\documentclass[12pt,oneside,reqno]{amsart}

\usepackage{amsfonts}
\usepackage{amsmath}
\usepackage{amssymb}
\usepackage{amsthm}
\usepackage{xcolor}
\usepackage{enumerate}
\usepackage{graphicx}
\usepackage{hyperref}
\usepackage{mathrsfs}
\usepackage{stmaryrd}
\usepackage[nameinlink,capitalize]{cleveref}

\newcommand{\be}{\begin{eqnarray}}
\newcommand{\ee}{\end{eqnarray}}
\newcommand{\ce}{\begin{eqnarray*}}
\newcommand{\de}{\end{eqnarray*}}

\newtheorem{theorem}{Theorem}[section]
\newtheorem{proposition}[theorem]{Proposition}
\newtheorem{lemma}[theorem]{Lemma}
\newtheorem{corollary}[theorem]{Corollary}
\newtheorem{remark}[theorem]{Remark}
\newtheorem{definition}[theorem]{Definition}

\newtheorem{assumption}[theorem]{Assumption}

\def\bt{\begin{theorem}}
\def\et{\end{theorem}}
\def\bp{\begin{proposition}}
\def\ep{\end{proposition}}
\def\bl{\begin{lemma}}
\def\el{\end{lemma}}
\def\bc{\begin{corollary}}
\def\ec{\end{corollary}}
\def\bd{\begin{definition}}
\def\ed{\end{definition}}
\def\br{\begin{remark}}
\def\er{\end{remark}}
\def\bx{\begin{Examples}}
\def\ex{\end{Examples}}
\def\ba{\begin{assumption}}
\def\ea{\end{assumption}}

\def\a{\alpha}
\def\b{\beta}

\def\[{{\Big[}}
\def\]{{\Big]}}
\def\<{{\langle}}
\def\>{{\rangle}}
\def\({{\Big(}}
\def\){{\Big)}}

\def\>{\geqslant}
\def\<{\leqslant}

\def\no{\nonumber}

\def\min{{\mathord{{\rm min}}}}
\def\max{{\mathord{{\rm max}}}}

\def\cB{{\mathcal B}}

\def\cF{{\mathcal F}}

\def\cH{{\mathcal H}}

\def\cW{{\mathcal W}}

\def\mB{{\mathbb B}}

\def\mE{{\mathbb E}}

\def\mM{{\mathbb M}}

\def\mP{{\mathbb P}}

\def\mR{{\mathbb R}}

\def\s{\sigma}
\def\Om{\Omega}
\def\om{\omega}
\def\la{\langle}
\def\ra{\rangle}
\def\hB{\hat{B}}
\def\hW{\hat{W}}
\def\hX{\hat{X}}
\def\hZ{\hat{Z}}
\allowdisplaybreaks

\newcommand{\st}[1]{{\color{red}#1}}
\newcommand{\sptext}[3]{\hspace{#1 em}\mbox{#2}\hspace{#3 em}}

\begin{document}

\title{Geiss-Besov Regularity of Exit Times of Diffusion Processes}

\author{Yao Chen$^{1}$, Jiagang Ren$^{1}$, Hua Zhang$^{2}$$^*$}

\dedicatory{$^{1}$School of Mathematics, Sun Yat-Sen University,\\
Guangzhou, Guangdong 510275, P.R.China\\
$^{2}$School of Statistics and Data Science, Jiangxi University of Finance and Economics,\\
Nanchang, Jiangxi 330013, P.R.China\\
Emails: Y. Chen: cheny2378@mail2.sysu.edu.cn\\
J. Ren: renjg@mail.sysu.edu.cn\\
H. Zhang: zh860801@163.com}


\thanks{This work is supported by National Natural Science Foundation of China (Grant Nos. 12261038 and 12371152), and Natural Science Foundation of Jiangxi Province (Grant Nos. 20232BAB201004, 20242BAB23003).}

\thanks{$^*$Corresponding Author.}

\begin{abstract}
We prove that the truncated exit times of possibly degenerate diffusion processes from open domains satisfying the exterior sphere condition belong to the Geiss-Besov space $\mB_p^{\frac{1}{2p}}$. In the one-dimensional case, the same result is proven using Zvonkin transformation for untruncated exit times if the diffusion is uniformly elliptic.
\end{abstract}

\maketitle
\def\ts{\tilde{\s}}
\def\tw{\tilde{w}}
\def\tB{\tilde{B}}
\def\tX{\tilde{X}}
\def\hs{\hat{\s}}
\section{Introduction}
A scale of Besov spaces on the classical Wiener space, which we call Geiss-Besov spaces below, has been introduced
in \cite{Geiss:Ylinen:21, Geiss26} to measure the fractional regularity of Wiener functionals. These spaces are
connected to, but different from, the traditional fractional Sobolev spaces on the abstract Wiener space which are defined by
S. Watanabe in \cite{wa} using
interpolation theory. They are motivated by the fact that many useful Wiener functionals have only fractional regularity; nevertheless, this
fractional regularity
is still very important in various aspects, e.g., in the theory and applications of backward
stochastic differential equations, see \cite{GGG:12,CGeiss:Steinicke:16, GT} and references therein. So it is a natural and important problem to determine the membership of typical Wiener functionals in
these spaces. This membership has been specified for local times of the one-dimensional  Brownian motion in \cite{Geiss:Ylinen:21}.
The purpose of the present paper is to address this problem for exit times of diffusion processes from open domains.
We would like to point out that the same problem
for non-degenerate diffusions in the context of Watanabe's fractional Sobolev spaces has been solved in \cite{amr}.

We organise the paper as follows. In Section 2 we state necessary background materials, including the definition of the Geiss-Besov spaces.
Our main results  are stated in Section 5. To prove these results, a key step is to estimate the moments of the difference of the exit times of the processes starting from different points and this will be done in Section 3 and Section 4, dealing with one- and multi-dimensional cases, respectively.

 \section{\texorpdfstring{Background Materials}{Background Material}}

Let $W_0:=C_0([0,\infty),\mR^m)$ be the space of continuous functions $w: [0,\infty)\mapsto\mR^m$ such that $w(0)=0$ and
let $\mu$ be the law of the standard $m$-dimensional Brownian motion on $W_0$. We equip $W_0$ with the compact uniform convergence
topology and denote by $\cW_0$ the completion of the associated Borel $\s$-algebra with respect to $\mu$. For $r\>0$, the shift operator
$\theta_r: W_0\mapsto W_0$ is defined by
$$
\theta_rw(\cdot):=w(\cdot+r)-w(r).
$$
Let $B$ and $B'$ be two independent $m$-dimensional Brownian motions  defined on some probability space
$(\Omega,\cF,\mP)$. Set for $0\<u<v<\infty$
$$ B^{u,v}(t)
   := \begin{cases}
      B(t)                             &: t\in [0,u] \\
      B(u) + \theta_uB'(t-u)        &: t\in (u,v] \\
      B^{u,v}(v)+ \theta_vB(t-v) &: t> v.
      \end{cases} $$
Then $B^{u,v}$ is again a Brownian motion. Consequently, if $\xi$ is a random variable on
$(W_0,\cW_0,\mu)$ then $\xi(B)$ and $\xi(B^{u,v})$ are well-defined random variables on $(\Omega,\cF,\mP)$
and the joint law of $(\xi(B), \xi(B^{u,v}))$ is independent of the choices of  $(\Omega,\cF,\mP)$,
$B$ and $B'$.  Hence the following definition, which is introduced in \cite{Geiss:Ylinen:21}, makes sense:

\bd[Geiss-Besov spaces] Let $\xi$ be a random variable on $(W_0,\cW_0,\mu)$ . For $p\>1$, $\alpha\in (0,1]$ define
$$
\|\xi\|_{p,\a}:=\sup_{0\<u<v<\infty}\frac{\|\xi(B)-\xi(B^{u,v})\|_p}{|u-v|^\a},
$$
and set
$$
\mB_p^\a:=\{\xi\in L_p:  \|\xi\|_{p,\a}<\infty\}.
$$
\ed

We refer to \cite{Geiss:Ylinen:21} and \cite{Geiss26} for a systematic study of these spaces.

Next we recall briefly the notion of strong solutions to stochastic differential equations, and the reader is referred to \cite{iw, rw} for more details.

Denote by $W$ the space of continuous functions from $\mR_+$ to $\mR^d$. Let $x_t$ be the canonical process on $W$
(resp. $W_0$) so that
$$
x_t(w):=w(t),~~t\>0,~w\in W~(\mbox{resp.} W_0).
$$
Denote by $\cH_t$ (resp. $\cH_t^0$) the $\s$-algebra on $W$ (resp. $W_0$) generated by $\{x_s, s\<t\}$  and $\{\bar{\cH}_t^0\}$ the
$\mu$-augmentation of $\{\cH^0_t\}$.

Given $m+1$ vector fields $b, \s_1,\cdots, \s_m$  on $\mR^d$,
a function
$$
F:\mR^d\times W_0\mapsto W
$$
is said to be a strong solution to the SDE $e(\s,b)$  if

(i)
$$
F^{-1}(\cH_t)\subset \cB(\mR^d)\times \bar{\cH}_t^0;
$$

(ii) for any $m$-dimensional $\{\cF_t\}$-Brownian motion $B=(B^1,\cdots, B^m)$
on any filtered probability space $(\Om, \cF, \{\cF_t\}, \mP)$ satisfying the usual conditions,  the process
$$
X(t,x):=F(x,B)(t)
$$
satisfies the following SDE:
\be\label{eq01}
\begin{cases}
&dX(t)=\s_i(X(t))dB^i(t)+b(X(t))dt,\\
&X(0)=x,
\end{cases}
\ee
where Einstein's convention for sums over repeated indices is used.

Now suppose the SDE $e(\s,b)$ has a unique strong solution. Let $D\subset\mR^d$ and $x\in\mR^d$. Define
$$
\tau_x(w):=\inf\{t\> 0: F(x, w)(t)\notin D\}.
$$
The main problem we will address is the following: for which $(p,\a)$ one has $\tau_x\in\mB_p^\a$ or $\tau_x^T\in\mB_p^\a$ where
$\tau_x^T:=\tau_x\wedge T$, and what is the best result one could expect to obtain? The answers will be given in Section 5.

Throughout the paper $C$ or $c$ with or without an index will denote a generic positive constant whose precise value is unimportant and
may vary from position to position.

\section{One-dimensional Diffusions}

In this section we will consider one-dimensional diffusions, so that $m=d=1$.

\subsection{One-dimensional Brownian Motion}

In this subsection we prepare some results for the one-dimensional Brownian motion. The aim is two-fold:
on the one hand, this is a good starting point for general cases; on the other hand, this is the case where one can make explicit
computations and thus gain insight into what is the optimal result we can expect to obtain.
\smallskip

We will need the following auxiliary lemma.

\bl\label{auxi}
Define for $u\in [0,1]$:
$$
\phi_p(u): = \frac{\Gamma(p+1)2^{p+2}}{\pi^{2p+1}} \sum_{k=0}^\infty\frac{\sin((2k+1)\pi u)}{(2k+1)^{2p+1}}.
$$
Then
there exist $c_p,C_p>0$ and $\delta \in (0,1)$ such that
$$
\phi_p(u) \< C_p u \sptext{.5}{for}{.5} u\in [0,1]
\sptext{1}{and}{1}
\phi_p(u)\> c_pu \sptext{.5}{for}{.5} u\in [0,\delta].
$$
\el

Now let $(\Om,\cF, \{\cF_t\}, \mP)$ be a probability space satisfying the usual conditions and $B$ be a one-dimensional standard
$\{\cF_t\}$-Brownian motion defined thereon.   Let $a<0<b$  and define
$$
\s^0_x(w):=\inf\{t\> 0: x+w(t)\in \{a,b\}\},
$$
$$
\s_x(\om):=\s_x^0(B(\om)).
$$

We have
\bigskip

\bt\label{th2.1}
For every $p\>1$ there exists a constant $C_p$ such that for all $x,y\in [a, b]$ one has
\be\label{eqn:th2.1}
    \mE[|\s_x-\s_y|^p] = (b-a)^{2p} \phi_p \left ( \frac{|x-y|}{b-a} \right )
    \< C_p|x-y|.
\ee
\et
\medskip

\begin{proof}
For the two cases $x=y$ and $[x,y]=[a,b]$ we have
$\s_x=\s_y$ and  $\phi_p \left ( \frac{y-x}{b-a} \right )=0$, so that we exclude
these two cases in the following. Suppose without any loss of generality that $x<y$.
For $z\in [a,b]$ and $w\in W_0$ set
$$
\tau_z^{0,1}(w):=\inf\{t\ge 0: z+w_t=a\},~~~\tau_z^{0,2}(w):=\inf\{t\ge 0: z+w_t=b\},
$$
$$
\tau_z^1(\om):=\tau^{0,1}_z(B(\om)),~~~\tau_z^2(\om):=\tau^{0,2}_z(B(\om)).
$$
We claim first that
\be\label{2}
|\s_x(\om)-\s_y(\om)|=\s^0_{a+y-x}(\theta_{\tau^1_x(\om)}B(\om))1_{\{\s_x<\s_y\}}+
\s^0_{b+x-y}(\theta_{\tau^2_y(\om)}B(\om))1_{\{\s_y<\s_x\}}.
\ee
If $\s_x(\om)<\s_y(\om)$, then $\tau_x^1(\om)<\tau_x^2(\om)$. Thus
$$
\s_x(\om)=\tau_x^1(\om)<\tau^2_y(\om).
$$
Consequently
$$
\s_y(\om)=\tau_x^1(\om)+\s_{a+y-x}^0(\theta_{\tau^1_x(\om)}B(\om)),
$$
from which we deduce that
\be\label{3}
0\<\s_y(\om)-\s_x(\om)=\s_{a+y-x}^0(\theta_{\tau^1_x(\om)}B(\om)).
\ee
If $\s_y(\om)<\s_x(\om)$, we have similarly
\be\label{4}
0\<\s_x(\om)-\s_y(\om)=\s_{b-y+x}^0(\theta_{\tau^2_y(\om)}B(\om)).
\ee
Combining (\ref{3}) and (\ref{4}) gives (\ref{2}).

Setting $u(t,x):=\mP(\s_x>t)$, we have (see e.g. \cite[Page 31]{im})
$$
u(t,x)=\frac{4}{\pi}\sum_{k=0}^\infty\frac{\exp(-\frac{(2k+1)^2\pi^2t}{2(b-a)^2})}{2k+1}
\sin \left ( \frac{(2k+1)\pi(x-a)}{b-a} \right ).
$$
This gives
$$
\mE[\s_x^p]=(b-a)^{2p} \phi_p\left (\frac{x-a}{b-a}\right ). $$
Since $\{\s_x<\s_y\}, \{\s_y<\s_x\}\in\cF_{\s_x\wedge\s_y}$, we have by the strong Markov property of Brownian motion
\ce
&   &  \mE[|\s_x-\s_y|^p] \\
& = &  \mE[|\s_{a+y-x}^0(\theta_{\tau^1_x}B)1_{\{\s_x<\s_y\}}|^p]
      +\mE[|\s_{b-y+x}^0(\theta_{\tau^2_y}B)1_{\{\s_y<\s_x\}}|^p] \\
& = & \mE[|\s_{a+y-x}^0(\theta_{\tau^1_x}B)|^p] \mP(\s_x<\s_y)
      +\mE[|\s_{b-y+x}^0(\theta_{\tau^2_y}B)|^p] \mP(\s_y<\s_x) \\
& = &  \mE[|\s_{a+y-x}|^p] \mP(\s_x<\s_y)
      +\mE[|\s_{b-y+x}|^p] \mP(\s_y<\s_x) \\
& = & (b-a)^{2p} \left [ \phi_p\left ( \frac{(a+y-x)-a}{b-a}\right ) \mP(\s_x<\s_y)
      +   \phi_p\left ( \frac{b-(b-y+x)}{b-a}\right )  \mP(\s_y<\s_x) \right ] \\
& = & (b-a)^{2p}  \phi_p\left ( \frac{y-x}{b-a}\right ) \left [ \mP(\s_x<\s_y)
      +   \mP(\s_y<\s_x) \right ].
\de
As our assumptions on $x$ and $y$ imply that $\mP(\s_y=\s_x) =0$ the assertion follows.
\end{proof}

\subsection{Non-degenerate Diffusion Case}

Now we turn to non-degenerate diffusions. We denote $\s_1$ simply by $\s$ and we make the following assumption:

\ba\label{ass:1d-SDE}
The functions
$\sigma,b$: $\mR\mapsto \mR$ are continuous. Moreover,  $b(\cdot)$ is Lipschitz continuous and compactly supported, and
there exists $c>0$ such that $|\s|\>c$ and there exists an increasing function $\rho: \mR_+\mapsto \mR_+$ such
that
$$
\int_{\mR_+}\rho^{-2}(u)du=\infty
\sptext{1}{and}{1}
|\s(x)-\s(y)|\<\rho(|x-y|).
$$
\ea

Consequently we have
$$
\int_{\mR}\s^{-2}(u)|b(u)|du<\infty.
$$

It is then a famous result of Yamada and Watanabe  \cite{yw} that Equation $e(\s, b)$  has a unique strong solution $F$.
Let $B$ be a one-dimensional Brownian motion on some filtered probability space $(\Om, \cF, \{\cF_t\}, \mP)$
satisfying the usual conditions and set
$$
X(t,x,\om):=F(x, B(\om))(t).
$$
We fix $-\infty < a < b < \infty$.  Define for $x\in [a,b]$
the stopping time
$$
\tau_x(\om):=\inf\{t \> 0: X(t,x,\om)\in \{a,b\}\}.
$$

\medskip

We have:
\bl\label{lemma3.1}
\begin{enumerate}[{\rm (1)}]
\item \label{item:2:lemma3.1}
      There exists a constant $C>0$ such that for any $\{\cF_t\}$-stopping time $\tau$ and $x,y\in\mR$ one has
      \be\label{lemma1dim1-1}
      \mE[|X(\tau,x)-X(\tau,y)|] \< C|x-y|. \ee
\item \label{item:3:lemma3.1}
      For any $p\>1$ there exists a constant $C>0$ such that
      \be \label{eqn:item:3:lemma3.1}
      \mE[\tau_x^p]\< C (|x-a|\wedge|x-b|).
      \ee
\end{enumerate}
\el
\bigskip

\begin{proof}
We first use the Zvonkin transformation (cf. \cite{z}). Set
$$
\psi(x):=\int_0^x\exp(-2\int_0^yb(u)\s^{-2}(u)du)dy.
$$
Then $\psi$ is a strictly increasing $C^2$-function on $\mR$ with
$$
0<\inf_{x\in\mR}\psi'(x)\<\sup_{x\in\mR}\psi'(x)<\infty.
$$
In particular, $\psi$ is bi-Lipschitz, i.e., there exists $c_\psi\> 1$ such that
$$
c_\psi^{-1}\<\frac{\psi(x)-\psi(y)}{x-y}\<c_\psi,~~~\forall x\neq y.
$$
Set again
$$
g(x):=(\s\psi')\circ \psi^{-1}(x).
$$
Taking $\gamma(x):=x$ in Lemma \ref{rhogamma} below,
one sees that there exists a function $\tilde{\rho}: [0,1]\mapsto \mR_+$ such that
$$
\int_{0+}\tilde{\rho}^{-2}(u)du=\infty,~~~|g(x)-g(y)|\<\tilde{\rho}(|x-y|).
$$
Consequently, by \cite{yw}, the equation
$$
Y(t,x)=\psi(x)+\int_0^t g(Y(s,x))dB_s
$$
has a unique solution $Y$. Then, by It\^o's formula,  the process $Z(t,x):=\psi^{-1}(Y(t,x))$ satisfies
$$
Z(t,x)=x+\int_0^t\s(Z(s,x))dB_s+\int_0^t b(Z(s,x))ds.
$$
Thus we have by the pathwise uniqueness that almost surely
$$
X(t,x)=\psi^{-1}(Y(t,x))~~\forall t\>0.
$$

\medskip

\eqref{item:2:lemma3.1} Set
$$
U_t:=Y(t,x)-Y(t,y),
$$
and denote by $L^0_t$  the local time of $U$ at $0$. Then we have by \cite{legall}
$$
L^0_t=0.
$$
Hence we have by Tanaka's formula
$$
|U_{\tau\wedge t}|=|\psi(x)-\psi(y)|+\int_0^{\tau\wedge t}\mbox{sgn}(U_s)(g(Y(s,x))-g(Y(s,y)))dB_s,
$$
which yields
$$
\mE[|U_{\tau\wedge t}|]=|\psi(x)-\psi(y)|.
$$
Letting $t\to\infty$ we have by Fatou's lemma
$$
\mE[|U_\tau|]\<|\psi(x)-\psi(y)|\<C|x-y|.
$$
This gives (\ref{lemma1dim1-1}) since
$$
c|Y(\tau,x)-Y(\tau,y)|\<|X(\tau,x)-X(\tau,y)|\<C|Y(\tau,x)-Y(\tau,y)|.
$$

\medskip

\eqref{item:3:lemma3.1} Obviously
$$
\tau_x=\inf\{t>0: Y(t,x)\in \{\psi(a),\psi(b)\}\}.
$$
Since there is a $c_g\ge 1$ such that
$$ \frac{1}{c_g} \< g^2(x) \< c_g
   \sptext{1}{for}{1} x\in \mR, $$
there exists a standard Brownian motion $\tB$ such that
$$
Y(t,x)=\psi(x)+\tB\left (\int_0^t g^2(Y(s,x))ds \right ).
$$
Setting
$$
\kappa_x:=\inf\{t>0: \psi(x)+\tB(t)\in \{\psi(a),\psi(b)\}\},
$$
we get $$
\tau_x\< c_g \kappa_x.
$$
Hence by Theorem \ref{th2.1} we may conclude with
\ce
     \mE[\tau_x^p]
&\<& C\mE[\kappa_x^p]\\
&\<& C |\psi(b)-\psi(a)|^{2p-1} (|\psi(x)-\psi(a)|\wedge|\psi(b)-\psi(x)|)\\
&\<& C|\psi(b)-\psi(a)|^{2p-1} (|x-a|\wedge|b-x|). \qedhere
\de
\end{proof}
\bigskip

Next we prove the lemma used in the proof.

\bl\label{rhogamma}
Let $\rho, \gamma: [0,1]\mapsto \mR_+$ be increasing and continuous and
$$
\rho(0)=\gamma(0)=0, ~~~\int_0^1\frac{d\gamma(x)}{\rho(x)}=\infty.
$$
Then
$$
\int_0^1\frac{d\gamma(x)}{\rho(x)+\gamma(x)}=\infty.
$$
\el
\begin{proof}
Take a decreasing sequence $\{x_n\}$ such that
$$
x_1=1,~~~\gamma(x_{n+1})=\frac12\gamma(x_n), ~~\forall n=1,2,\cdots.
$$
Then
$$
\frac12\sum_{n=1}^\infty\rho^{-1}(x_n)\gamma(x_n)
\<\int_0^1\frac{d\gamma(x)}{\rho(x)}\<\sum_{n=1}^\infty\rho^{-1}(x_n)\gamma(x_n).
$$
Hence
$$
\int_0^1\frac{d\gamma(x)}{\rho(x)}=\infty\Longleftrightarrow
\sum_{n=1}^\infty\rho^{-1}(x_n)\gamma(x_n)=\infty.
$$
Set
$$
\mM:=\{n: \rho(x_n)\<\gamma(x_n)\}.
$$
If $\sharp(\mM)=\infty$, then
$$
\int_0^1\frac{d\gamma(x)}{\rho(x)+\gamma(x)}\>
\frac14\sum_{n\in\mM}\frac{\gamma(x_n)}{\gamma(x_n)}=\infty.
$$
If $\sharp(\mM)<\infty$, then there exists an $N$ such that
$$
\int_0^1\frac{d\gamma(x)}{\rho(x)+\gamma(x)}\>
\frac14\sum_{n\>N}\frac{\gamma(x_n)}{\rho(x_n)}=\infty.
$$

This completes the proof.
\end{proof}

\bigskip

We are now ready to prove the main result in this subsection.

\bt\label{main result one dimensional homogeneous}
Let $p\>1$ and $-\infty<a<b<\infty$. Then there exists a constant $C$ such that for any $x,y\in (a,b)$
 one has
$$
\mE[|\tau_x-\tau_y|^p]\<  C|x-y|
$$
\et

\begin{proof}
	Suppose without any loss of generality
	that $x<y$.
	For $z\in [a,b]$ and $w\in W_0$ set
	\ce
	\tau_z^0(w)&:=&\inf \{ t\> 0; F(z,w)(t) \in \{a,b\}\},\\
	\tau_z(\omega)&:=&\tau_z^0(B(\omega)).
	\de
	We have
	\ce
	|\tau_x-\tau_y|
	=
	\tau^0_{X(\tau_x,y)}
	(\theta_{\tau_x}B(\omega))
	1_{\{\tau_x<\tau_y\}}
	+
	\tau^0_{X(\tau_y,x)}
	(\theta_{\tau_y}B(\omega))
	1_{\{\tau_y<\tau_x\}}.
	\de
Since $\{\tau_x< \tau_y \},\{\tau_y< \tau_x \} \in \cF_{\tau_x \wedge \tau_y}$,  by the strong Markov property and Lemma \ref{lemma3.1}
we have
	\ce
\mathbb E|\tau_x-\tau_y|^p
&=&
\mathbb E\left[
\left(
\tau^0_{X(\tau_x,y)}
(\theta_{\tau_x}B)
\right)^p
1_{\{\tau_x<\tau_y\}}
\right]
+
\mathbb E\left[
\left(
\tau^0_{X(\tau_y,x)}
(\theta_{\tau_y}B)
\right)^p
1_{\{\tau_y<\tau_x\}}
\right]\\
	&= &
	\mathbb E\left[
	1_{\{\tau_x<\tau_y\}}
	\left.
	\mathbb E\left[
	\left(\tau^0_z(B)\right)^p
	\right]
	\right|_{z=X(\tau_x,y)}
	\right]+\mathbb E\left[
	1_{\{\tau_y<\tau_x\}}
	\left.
	\mathbb E\left[
	\left(\tau^0_z(B)\right)^p
	\right]
	\right|_{z=X(\tau_y,x)}
	\right]\\
	&\<&
	C\mathbb E\left[
	(X(\tau_x,y)-a)
	1_{\{\tau_x<\tau_y\}}
	\right]+	C\mathbb E\left[
	(b-X(\tau_y,x))
	1_{\{\tau_y<\tau_x\}}
	\right].
	\de
	By \cite[Theorem 1.1]{iw},
	if $x<y$, we have almost surely
	\ce
	X(t,x)\< X(t,y),
	\qquad t\ge0.
	\de
	Hence if  $\{\tau_x<\tau_y\}$,
	then $ X(\tau_x,x)=a$, and if  $\{\tau_y<\tau_x\}$, then $ X(\tau_y,y)=b$.
	Combining the preceding estimates and Lemma \ref{lemma3.1},
	we obtain
	\ce
	\mathbb E[|\tau_x-\tau_y|^p]
	&\< &
	C\mathbb E\left[
	|X(\tau_x,y)-X(\tau_x,x)|
	1_{\{\tau_x<\tau_y\}}
	\right]\\
	&&+
	C\mathbb E\left[
	|X(\tau_y,y)-X(\tau_y,x)|
	1_{\{\tau_y<\tau_x\}}
	\right]\\
	&\<&
	C\mathbb E\left[
	|X(\tau_x\wedge\tau_y,y)
	-X(\tau_x\wedge\tau_y,x)|
	\right]\\
	&\<&
	C|x-y|.
	\de
\end{proof}

\section{Multidimensional Diffusion: Truncated Exit Time}
\label{sec:Rd-case}

We now turn to the multidimensional case. For $\s=(\s_1,\cdots, \s_m)$ and  $b$ we consider the Equation $e(\s, b)$ and we make the following assumption throughout this section:

\ba\label{assumption:rd}
\begin{enumerate}[{\rm (1)}]
\item Suppose $\s_i, b: \mR^d\mapsto\mR^d$, $i=1,\cdots,m$ are $m+1$ bounded continuous vector
      fields on $\mR^d$, and assume
$$
K:=\sup\left \{\sum_{i=1}^m\frac{|\s_i(x)-\s_i(y)|^2}{|x-y|^2}+2\frac{|b(x)-b(y)|}{|x-y|}:
x\neq y\right \}<\infty.
$$

\item

$D$ is a bounded open domain satisfying a uniform exterior sphere condition, i.e., there exists an $R>0$ such that for every
$y\in\partial D$ there is a $z\in \mR^d$ such that $\overline{B(z, R)}\cap\overline{D}=\{y\}$.

\item
There are $\delta,c>0$ such that for all $x\in B(\delta)\cap D$ there is a $y\in \pi(x)$  where
$$
\pi(x):=\{y\in\partial D: |x-y|=d(x,\partial D)\},
$$
such that
$$  \sum_{i=1}^m\left | \left (\s_i(x'), \frac{(y-x)}{|y-x|}\right) \right |^2
\> c^2 \sptext{1}{for}{1} x'\in \overline{D}\cap \overline{B(y,2 \delta)}.$$
\end{enumerate}
\ea
\bigskip

It is well known that under this assumption the Equation $e(\s,b)$ admits a unique strong solution $F(x,w)$. Let $B=(B^1,\cdots, B^m)$ be any $m$-dimensional Brownian motion
on any filtered probability space satisfying the usual conditions. Then the process
$$
X(t,x,\om)=F(x,B(\om))(t)
$$
satisfies the equation $(\ref{eq01})$ and the following result is well known, see, e.g., \cite{iwb}.

\bigskip

\bl\label{statement:sensitivity:SDE:starting_value}
$X(t,x)$  constitutes a stochastic flow of homeomorphisms, and  there exist two positive constants $C$ and $K$ such that
$$
\mE[|X(\tau,x)-X(\tau,y)|]\< C|x-y|e^{KT}
$$
for all $T\> 0$ and all stopping times $\tau:\Omega \to [0,T]$.
\el
\medskip

The proof of the following lemma is standard and therefore omitted:
\medskip

\bl\label{statement:straight_line}
For every $x\in D$ and every $y\in \pi(x)$,
there is only one $z\in \mR^d$ such that
$$
\overline{B(z, R)}\cap\overline{D}=\{y\}.
$$
Moreover, $x,y,z$ belong to the same straight line.
\el
\bigskip

We will also need the following lemma.

\bigskip
\bl\label{lemma41_new}
Let $ \xi=(\xi(t))_{t\> 0}$, $\xi(t)=(\xi_1(t),\ldots,\xi_d(t)):\Omega \to \mR^d$, be such that
$$ \xi_j(t) :=M_j(t)+\int_0^t\beta_j(s)ds, $$
where $(M_j(t))_{t\> 0}$ are continuous martingales such that $M_j(0)\equiv 0$.
Let $-b < 0 < a$ and $h>0$, and
$Q:= [-b,a] \times [-h,h]^{d-1}$. Assume that
$[M_1 ]_{t\> 0}$ is strictly increasing, and there is a constant $K_1\> 1$ such that
$$ K^{-1}_1 t\<[M_1]_t\<K_1t
   \sptext{1}{and}{1}
   [ M_j ]_t\<K_1t \sptext{1}{for}{1} t\> 0 \sptext{.5}{and}{.5} j=2,\ldots,d $$
and $(\beta_j(t))_{t\ge 0}$ are progressively measurable processes such that
$$
\sup_{t\>0}\{|\beta_j(t)|\}\<K_2
$$
for some constant $K_2>0$. Let
\begin{align*}
\sigma^{O,T}  & :=\inf\{t\>0: \xi_1(t) = a\}\wedge T,\\
\sigma^{I,T}_1& :=\inf\{t\>0: \xi_1(t) = -b\}\wedge T, \\
\sigma^{I,T}_j& :=\inf\{t\>0: |\xi_j(t)|=h\}\wedge T \sptext{1}{for}{1} j=2,\ldots,d,\\
\sigma^{I,T} &:= \sigma^{I,T}_1 \wedge \sigma^{I,T}_2 \wedge \cdots \wedge \sigma^{I,T}_d.
\end{align*}
Then for every $T>0$ and $p\>1$ there is a $C>0$ such that
one has
\be\label{eqn:lemma41}
     \mP(\sigma^{I,T} < \sigma^{O,T})
  \< C a
\ee
and
\be
    \mE [  |\sigma^{O,T}|^p]
 \< Ca.
\ee
\el

\begin{proof}
(a) Let $B$ be a Brownian motion such that $M_1(t)=B([M_1]_t)$. Then
$$
\xi_1(t)=B(\la M_1\ra_t)+\int_0^t\beta_1(s)ds.
$$
For $t<T$ we have
\ce
     \mP(\sigma^{O,T}\<t)
& =& \mP(\max_{0\<s\<t}\xi_1(s)\>a)\\
&\>& \mP(\max_{0\<s\<t}(B(\la M_1 \ra_s)-K_2s)\>a)\\
&\>& \mP(\max_{0\<s\<t}(B(K_1^{-1}s)-K_2s)\>a)\\
& =& \mP(\max_{0\<s\<K_1^{-1}t}(B_s-K_2K_1s)\>a).
\de
Since the last quantity is equal to (see, e.g., \cite[p. 250]{borodin})
$$
\Phi\left (\frac{a}{\sqrt{2K_1^{-1}t}}+\frac{K_2\sqrt{K_1t}}{\sqrt{2}}\right )
+e^{-2K_1K_2a}\Phi\left (\frac{a}{\sqrt{2K_1^{-1}t}}-\frac{K_2\sqrt{K_1t}}{\sqrt{2}}\right ),
$$
where
$$
\Phi(x):=\frac{1}{\sqrt{\pi}}\int_x^\infty e^{-u^2}du,
$$
noticing $\Phi(-\infty)=1$ we have
$ 1-\Phi(x)=\Phi(-x)$.
Hence
\be\label{taudayut}
     \mP(\sigma^{O,T}>t)
&\<& \left[ \Phi\left (-\frac{a}{\sqrt{2K_1^{-1}t}}-\frac{K_2\sqrt{K_1t}}{\sqrt{2}}\right)
           -\Phi\left (\frac{a}{\sqrt{2K_1^{-1}t}}-\frac{K_2\sqrt{K_1t}}{\sqrt{2}} \right) \right]\no\\
&  & +(1-e^{-2K_1K_2a})\Phi\left (\frac{a}{\sqrt{2K_1^{-1}t}}-\frac{K_2\sqrt{K_1t}}{\sqrt{2}}\right )\no\\
&\<&c(t^{-\frac12}+1)a
\ee
with $c=c(T,K_1,K_2)>0$. Thus
$$
\mE[\sigma^{O,T}] \< \int_0^{T} \mP(\sigma^{O,T}>t)dt
\< c \int_0^T (t^{-\frac{1}{2}} +1) dt \, a
 = c\, (2\sqrt{T} + T) \, a.
$$
Hence
\ce
     \mP(\sigma^{O,T}>\sigma_1^I)
&\<& \mP\left (\inf_{0\<t\<\sigma^{O,T}}\xi_1(t)=-b\right )
 \<  \mP\left (\max_{0\<t\<\sigma^{O,T}}\left |M_1(t)+\int_0^t\beta_1(s)ds\right |\>b\right)\\
&\<& \mP\left (\max_{0\<t\<\sigma^{O,T}}|M_1(t)|\>\frac{b}{2}\right )+\mP\left(\max_{0\<t\<\sigma^{O,T}}\left |\int_0^t\beta_1(s)ds\right |\>\frac{b}{2}\right )\\
&\<& \frac{4}{b^2} \mE \[ \max_{0\<t\<\sigma^{O,T}}|M_1(t)|^2\] +\mP\left(\sigma^{O,T} K_2 \> \frac{b}{2} \right ) \\
&\<& \frac{16}{b^2} \mE \[ \langle M_1 \rangle_{\sigma^{O,T}}\] + \frac{2K_2}{b} \mE[  \sigma^{O,T}] \\
&\<& \left [ \frac{16}{b^2} K_1 + \frac{2K_2}{b} \right ] \mE [\sigma^{O,T}] \\
&\<& \left [ \frac{16}{b^2} K_1 + \frac{2K_2}{b} \right ] c\, (2\sqrt{T} + T) \, a.
\de
Similarly, for $j=2,\ldots,d$ we derive
$$
    \mP(\sigma^{O,T}>\sigma_j^I)
 \< \mP\left (\max_{0\<t\<\sigma^{O,T}}\left |M_j(t)+\int_0^t\beta_j(s)ds\right |\>h\right),
$$
so that the same argument applies.
\smallskip

(b) We have by \eqref{taudayut} that
$$
   \mE[|\sigma^{O,T}|^p] = p\int_0^Tt^{p-1}\mP(\sigma^{O,T}>t)dt
\< \left [ p c \int_0^T(t^{-\frac12}+1)t^{p-1}dt \right ] \, a. \qedhere
$$
\end{proof}

\medskip

Set
$$
I(\delta) := \{ x \in D : d(x,\partial D) > \delta \}
\sptext{1}{and}{1}
B(\delta) := \{ x \in \overline{D} : d(x,\partial D) \< \delta \}.
$$
\medskip

For an arbitrarily fixed $x\in B(\delta)$,
up to an orthogonal transformation we can assume that
$$ x=(-\xi,0,\ldots,0), \quad
   y = 0, \quad
   z =(R,0,\ldots,0), \quad
   \gamma:=(1,0,\ldots,0). $$
Let $a,b\in (0,R)$ such that $a^2+b^2=R^2$. We define
$u:[-\delta,a]\times \overline{B(0,b;\mR^{d-1})} \to \mR^d$ by
$$ (x_1,\ldots,x_d) \mapsto \left ( \delta \left ( \frac{\delta+x_1}{\delta  + R -\sqrt{R^2-x_2^2-\cdots-x_d^2}}-1\right ),x_2,\ldots,x_d \right ),$$
where $B(0,b;\mR^{d-1})$ is the open centered ball with radius $b$ in $\mR^{d-1}$.
\medskip

Let $$ 0< \varepsilon \< \frac{1}{2} \frac{c}{\sqrt{m} \sup_{i=1,\ldots,m} \|\sigma_i\|_\infty} $$
and choose $r\in (0,\delta\wedge b)$ such that for $x'\in B(0,r)$ one has
$$ |\nabla u^1(x')-\gamma| < \varepsilon, $$
which implies
\be\label{positivity}
  &&   \left ( \sum_{i=1}^m \left | \left ( \nabla u^1(x'),\sigma_i(x') \right )\right |^2 \right )^\frac{1}{2}\no\\
&\>&  \left ( \sum_{i=1}^m \left | \left ( \gamma,\sigma_i(x') \right ) \right |^2 \right )^\frac{1}{2}
   - \left ( \sum_{i=1}^m \left |\left ( \nabla u^1(x')-\gamma,\sigma_i(x') \right ) \right |^2 \right )^\frac{1}{2}\no \\
&\>& c - \varepsilon \left ( \sum_{i=1}^m |\sigma_i(x')|^2 \right )^\frac{1}{2}
 \> \frac{c}{2}.
\ee

\medskip

For convenience we
extend
$u$ to $u=(u^1,\ldots,u^d):\mR^d\to \mR^d$ such that $u^1,\ldots,u^d\in C_b^2(\mR^d)$,
and define
$$Z(t) := u(X(t,x)).$$
Then for $j=1,\ldots,d$ one has by It\^o's formula
\ce
    Z_j(t)
  &=& u^j(X(t,x))\\
& =& u^j(x) + \int_0^t \left(\nabla u^j(X(s,x))\sigma(X(s,x))\right)_i d B^i(s)
    + \int_0^t \nabla u^j(X(s,x))b(X(s,x)) ds\\
& &\quad + \frac{1}{2} \sum_{i=1}^m \int_0^t \sigma_i^\intercal(X(s,x)) (Hu^j)(X(s,x)) \sigma_i(X(s,x)) ds,
\de
where $H u^j$ is the Hessian matrix of $u^j$, and we may use
$$ [Z_j]_t = \int_0^t \sum_{i=1}^m \left |\left ( \nabla u^j(X(s,x)),\sigma_i(X(s,x)) \right ) \right |^2 ds$$
as a version of the quadratic variation.

Now fix an $h>0$ such that
$$ Q := [-2h,0] \times [-h,h]^{d-1} \subseteq u(B(0,r))$$
and set
$$
\delta':= \min \{\delta, |u^{-1}((-h,0,\ldots,0))| \}.
$$

Using (\ref{positivity}) we have, by construction, that
$$ [Z_1]_{t\wedge \rho} = \int_0^{t\wedge \rho} \sum_{i=1}^m \left |\left ( \nabla u^1(X(s,x)),\sigma_i(X(s,x)) \right ) \right |^2 ds
   \> (t\wedge\rho) \frac{c^2}{4} $$
for any stopping time $\rho\< \sigma$ with
$$ \sigma:= \inf \{ t\ge 0: Z(t)\in \partial Q\}.$$

\medskip

For $x=(-\xi,0,\ldots,0)$ with $\xi \in (0,\delta']$\st{,} define
$$ F(x) := u^{-1} \Big ([-{2h},0] \times [-h,h]^{d-1}\Big ) \cap \overline{D},$$
$$
a:=|u((-\xi,0,\ldots,0))|,
$$
$$
S^O:=\partial F(x)\cap D^c\sptext{1}{and}{1}S^I:=\partial F(x)\cap D.
$$
For $T>0$ we
define the hitting times
\begin{align*}
\sigma^{O,T}_x&:= \inf \{ t\> 0:  X(t,x) \in S^O \}\wedge T,\\
\sigma^{I,T}_x&:= \inf \{ t\> 0:  X(t,x) \in S^I \}\wedge T
\end{align*}
with $\sigma^{I,T}_x:=T$ if $S^I= \emptyset$.

\medskip

We have
\bl\label{statement:expeted_exit_time_d-dim-case}
Let $p\in [1,\infty)$. Then there exists a constant $C> 0$ such  that
$$ \mP(\sigma_x^{I,T} < \sigma_x^{O,T})  \< Cd(x,\partial D)
\sptext{1}{and}{1}
\mE  \[ |\sigma_x^{O,T}|^p 1_{\{\sigma_x^{O,T}\<\sigma_x^{I,T}\}} \] \< Cd(x,\partial D).
$$
\el
\begin{proof}
We let
$b:= 2h-a\>h$.
We redefine the diffusion coefficient of $Z_1$ to be
$$ \sigma_1^{\tilde Z}(t,\omega)
   := \begin{cases}
      \nabla u^1 (X(t,\omega)) \sigma(X(t,\omega)) &: Z(t,\omega)\in Q \\
      (1,0,\ldots,0)                               &: Z(t,\omega)\not \in Q
      \end{cases}.$$
The new process is denoted by $(\tilde Z(t))_{t\ge 0}$.
\medskip
Set
$$ \xi(t)     := \tilde Z(t) - \tilde Z(0) $$
and define
\begin{align*}
\gamma^O &: = \inf \{ t\ge 0: \tilde Z_1(t) = 0 \}\wedge T,\\
\gamma_1^I & := \inf \{ t\ge 0: \tilde Z_1(t) = -2h \}\wedge T, \\
\gamma_j^I & := \inf \{ t\ge 0: |\tilde Z_j(t)| = h \}\wedge T \sptext{1}{for}{1} j=2,\ldots,d,\\
\gamma^I &:=\gamma_1^I\wedge\cdots\wedge \gamma_d^I.
\end{align*}
By \cref{lemma41_new} we get
\begin{align*}
    \mE [|\gamma^O|^p]
&\< C|a|,\\
    \mP\Big (\gamma^{I} < \gamma^{O}\Big )
&\< C|a|.
\end{align*}
Consequently
$$\mP(\s_x^{I,T} < \s_x^{O,T})
\< \mP(\gamma^{I} < \gamma^{O})
\< C|a| $$
and
$$ \mE \[|\s_x^{O,T}|^p 1_{\{ \s_x^{O,T} \< \s_x^{I,T}\}}\]
  \<  \mE [|\gamma^{O}|^p ]
  \< C|a|.
$$
Since $u^1(0)=0$ and $u^1$ is Lipschitz, we have $|a|\<C|\xi|\<Cd(x,\partial D)$  and the desired estimates follow.
\end{proof}

\bc
\label{statement:general_upper_bound_EtauxTp}
For $x\in B(\delta)$ one has
\be \label{eqn:statement:general_upper_bound_EtauxTp}
   \mE  [|\tau_x^T|^p] \< Cd (x,\partial D).
\ee
\ec
\medskip

\begin{proof}
It suffices to notice that
$$
\tau_x^T\<\s_x^{O,T}1_{\{\s_x^{O,T}\<\s_x^{I,T}\}}+T1_{\{\s_x^{I,T}<\s_x^{O,T}\}}
$$
and use the above lemma.
\end{proof}

\bl \label{statement:general_estimate_taux-tauy}
For $\delta,T>0$, $p\in [1,\infty)$, and $x,y\in \overline{D}$ one has
\begin{align*}
     \mE \[ |\tau_x^T-\tau_y^T|^p \]
& \< 2^p \Bigg [ \[ \mE |\tau_z^T|^p 1_{\left \{\tau_x \< \tau_y^T, z= X(\tau_x,y) \in B(\delta)\right \}} \]
                +  \mE \[|\tau_z^T|^p 1_{\left \{\tau_y \< \tau_x^T, z= X(\tau_y,x) \in B(\delta)\right \}}\] \\
& \hspace{2.7em} +
T^p \mP\left (\tau_x \< \tau_y^T, X(\tau_x,y) \in I(\delta) \right ) + T^p \mP\left (\tau_y \< \tau_x^T,  X(\tau_y,x) \in I(\delta) \right )
\Bigg ].
\end{align*}
\el

\begin{proof}
We observe that
\ce
|\tau_x^T-\tau_y^T|^p
& \<& |\tau_x -\tau_y^T|^p 1_{\{\tau_x \< T \< \tau_y\}}
   +  |\tau_x^T -\tau_y|^p 1_{\{\tau_y \< T \< \tau_x\}} \\
&& +  |\tau_x -\tau_y|^p 1_{\{\tau_x \< \tau_y \< T\}}
   +  |\tau_x -\tau_y|^p 1_{\{\tau_y \< \tau_x \< T\}}
\de
and
$$
     |\tau_x -\tau_y| 1_{\{\tau_x \< \tau_y \< T\}}
  =  |\tau_x -\tau_y^T| 1_{\{\tau_x \< \tau_y \< T\}}. $$
This implies
$$
|\tau_x^T-\tau_y^T|^p
\< 2^p \left [ |\tau_x -\tau_y^T|^p 1_{\{ \tau_x \< \tau_y^T\}}
             + |\tau_y -\tau_x^T|^p 1_{\{ \tau_y \< \tau_x^T\}} \right ].
$$
Then
$$
  \mE \[|\tau_x -\tau_y^T|^p 1_{\{ \tau_x \< \tau_y^T\} \cap \{ X(\tau_x,y) \in B(\delta)\} }\]
 = \mE \[|\tau_z^T|^p 1_{\left \{\tau_x \< \tau_y^T, z= X(\tau_x,y) \in B(\delta)\right \}}\]
$$
and
$$
    \mE \[ |\tau_x -\tau_y^T|^p 1_{\{ \tau_x \< \tau_y^T\} \cap \{ X(\tau_x,y) \in I(\delta)\} }\]
\< T^p \mP\left (\tau_x \< \tau_y^T, X(\tau_x,y) \in I(\delta)\right ). \qedhere
$$
\end{proof}
\medskip

\bl\label{statement:conditional_upper_bound_EtauxTp}
For $x,y\in  \overline{D}$
one has
\be  \label{eqn:statement:conditional_upper_bound_EtauxTp}
\mE \left[|\tau_z^T|^p 1_{\left \{\tau_x \< \tau_y^T, z= X(\tau_x,y) \in B(\delta)\right \}}\right]
   \< C
       \mE \left[| X(\tau_x,y) - X(\tau_x,x)| 1_{\{ \tau_x \< \tau_y^T, X(\tau_x,y) \in B(\delta)\}}\right].
\ee
\el

\begin{proof}
We get that
\begin{align*}
     \mE \left[|\tau_z^T|^p 1_{\left \{ \tau_x \< \tau_y^T, z= X(\tau_x,y) \in B(\delta)\right \}}\right]
&\< C\mE  \left [ d (z,\partial D) 1_{\left \{\tau_x \< \tau_y^T, z= X(\tau_x,y) \in B(\delta)\right \}} \right ] \\
&\< C \mE  \left [ |X(\tau_x,y)-X(\tau_x,x)|
1_{\left \{\tau_x \< \tau_y^T, X(\tau_x,y) \in B(\delta)\right \}} \right ].\qedhere
\end{align*}
\end{proof}

Summarizing the above, we get the following estimate:

\bt\label{statement:upper_general_bound_tauxT-tauyT}
For $x,y\in \overline{D}$ one has
\be\label{eqn:statement:upper_general_bound_tauxT-tauyT}    \mE \left[|\tau_x^T-\tau_y^T|^p\right]
  \< C
   \Big [ \mE |X(\tau_x^T\wedge \tau_y^T,x)-X(\tau_x^T\wedge \tau_y^T,y)| \Big ].
\ee
\et

\begin{proof}
From Lemmas \ref{statement:general_estimate_taux-tauy} and
\ref{statement:conditional_upper_bound_EtauxTp}
we derive
\begin{align*}
     \mE \left[|\tau_x^T-\tau_y^T|^p\right]
& \< C \Bigg [  \mE\[ |\tau_z^T|^p 1_{\left \{\tau_x \< \tau_y^T,
          z= X(\tau_x,y) \in B(\delta)\right \}}\] \\
& \hspace{2.7em} +  \mE\[ |\tau_z^T|^p 1_{\left \{\tau_y \< \tau_x^T, z= X(\tau_y,x) \in B(\delta)\right \}} \]\\
& \hspace{2.7em} +
T^p \mP\left (| X(\tau_x,x) - X(\tau_x,y)| > \delta,
    \tau_x \< \tau_y^T, X(\tau_x,y) \in I(\delta) \right ) \\
& \hspace{2.7em} + T^p \mP\left (| X(\tau_y,y) - X(\tau_y,x)| > \delta,
    \tau_y \< \tau_x^T,  X(\tau_y,x) \in I(\delta) \right )
\Bigg ] \\
& \< C\Bigg [
       \mE \[| X(\tau_x,y) - X(\tau_x,x)| 1_{\{ \tau_x \< \tau_y^T, X(\tau_x,y) \in B(\delta)\}}\]\\
& \hspace{2.7em} +
       \mE \[ | X(\tau_y,x) - X(\tau_y,y)| 1_{\{ \tau_y \< \tau_x^T, X(\tau_y,x) \in B(\delta)\}} \]\\
& \hspace{2.7em} +
T^p \mP\left (| X(\tau_x,x) - X(\tau_x,y)| > \delta,
    \tau_x \< \tau_y^T, X(\tau_x,y) \in I(\delta) \right ) \\
& \hspace{2.7em} + T^p \mP\left (| X(\tau_y,y) - X(\tau_y,x)| > \delta,
    \tau_y \< \tau_x^T,  X(\tau_y,x) \in I(\delta) \right )
\Bigg ] \\
& \< C\Bigg [
       \mE \[| X(\tau_x^T\wedge \tau_y^T,y) - X(\tau_x^T\wedge \tau_y^T,x)|\]\\
& \hspace{3em} +
T^p \mP\left (| X(\tau_x^T\wedge \tau_y^T,x) -
                X(\tau_x^T\wedge \tau_y^T,y)| > \delta \right )
\Bigg ]\\
&\< C\mE \[| X(\tau_x^T\wedge \tau_y^T,y) - X(\tau_x^T\wedge \tau_y^T,x)|\].
\end{align*}
\end{proof}

We are now ready to state our main result in this section.
\bigskip

\bt\label{main result multidimensional truncated}
For every $T>0$ and $p\>1$ there exists a constant $ C= C(T,p)>0$ such that
$$
\mE[|\tau_x^T-\tau_y^T|^p]\<  C |x-y| \sptext{1}{for all}{1} x,y\in \overline{D}.
$$
\et

\begin{proof}
The result follows from
\cref{statement:upper_general_bound_tauxT-tauyT},
\cref{statement:expeted_exit_time_d-dim-case}, and
\cref{statement:sensitivity:SDE:starting_value}.
\end{proof}

\begin{remark}
In \cite[Lemma 9]{buckdahn} the authors consider situations when $D$ satisfies the exterior sphere condition and $X$ is
a uniformly elliptic diffusion, for which there is a $\theta<0$ such that
$$
\sup_{x,y\in \overline{D}\atop x\not= y} \frac{\mE [|e^{\frac{\theta}{2} \tau_x} - e^{\frac{\theta}{2} \tau_y}|]}{|x-y|}< \infty.
$$
\end{remark}

\section{Geiss-Besov Regularity}

Using the results obtained above, we can prove the Besov regularity of exit times.

We begin with the untruncated exit time in the one-dimensional case and suppose we are in the situation of Subsection 3.2. So let $m=d=1$ and let $F(x,w)$ be the unique strong solution of Equation $e(\s,b)$ where
$\s$ and $b$ satisfy Assumption \ref{ass:1d-SDE}. Fix $a<b$ and for $x\in [a,b]$ define
$$
\s^0_x(w):=\inf\{t\> 0: F(x,w)(t)\in \{a,b\}\}.
$$
Let $B$ and $B'$ be two independent  Brownian motions defined on some filtered probability space $(\Om, \cF, \{\cF_t\}, \mP)$
satisfying the usual conditions. For a fixed $r>0$ we set
$$
\hB_t=\begin{cases}
B'_t, &~0\< t\< r,\\
B'_r+\theta_rB_{t-r},&~t>r.
\end{cases}
$$
Set again
$$
X(t,x,\om):=F(x,B(\om))(t),~~~\hX(t,x,\om):=F(x,\hB(\om))(t).
$$
Then
\be
X(t,x)=x+\int_0^t\s(X(s,x))dB(s)+\int_0^tb(X(s,x))ds,
\ee
and
\be
\hX(t,x)=x+\int_0^t\s(\hX(s,x))d\hB(s)+\int_0^tb(\hX(s,x))ds.
\ee
Define
$$
\s_x:=\inf\{t\> 0: X(t,x)\in\{a,b\}\},~~~\hs_x:=\inf\{t\> 0:\hX(t,x)\in\{a,b\}\}.
$$
Then
$$
\s_x(\om)=\s^0_x(B(\om)),~~~\hs_x(\om)=\s_x^0(\hB(\om)).
$$
We now prove
\bl\label{lemma61}
For every $p\>1$ there exists a constant $C=C(p)$ which is independent of $r$ and $x$ such that
\be\label{besov1dim0}
\mE[|\s_x-\hs_x|^p]\<Cr^{\frac12},~~~~~~\forall x\in (a,b), ~~r\in [0,1].
\ee
\el
\begin{proof}
Set $\tB:=\theta_rB$ and
$$
\tX(t,x,\om):=F(x,\tB(\om))(t).
$$
Define
$$
\ts_x:=\inf\{t  \>0: \tX(t,x)\in\{a,b\}\}.
$$
Then
$$
\ts_x(\om)=\s_x^0(\tB(\om)).
$$
We use the following decomposition:
\be\label{besov1dim1}
|\s_x-\hs_x|&\<&|\ts_{X(r,x)}-\ts_{\hX(r,x)}|1_{\{\s_x\wedge\hs_x>r\}}+|\s_x|1_{\{\hs_x\<r, \s_x>r\}}\no\\
&&+|\hs_x|1_{\{\hs_x>r, \s_x\<r\}}+r1_{\{\s_x\vee\hs_x\<r\}}\\
&=:&\sum_{i=1}^4I_i.
\ee
We first have by the strong Markov property
\be\label{besov1dim2}
\mE[|I_1|^p]\<\mE\left[\left.\mE[|\ts_y-\ts_z|^p]\right|_{y=X(r,x), z=\hX(r,x)}\right]\<C\mE[|X(r,x)-\hX(r,x)|]\<Cr^{\frac12}.
\ee
Secondly, it is trivial that
\be\label{besov1dim3}
\mE[|I_4|^p]\<Cr^p.
\ee
It remains to treat $I_2$ and $I_3$. For $I_2$ we have by (\ref{eqn:item:3:lemma3.1})
\ce
\mE[|I_2|^p]\<\mE[|\s_x|^p]\mP(\hs_x\<r)\<C(|x-a|\wedge|x-b|)P(\hs_x\<r).
\de
But
\ce
\mP(\hs_x\<r)\<\mP(\max_{0\<s\<r}|\hX(s,x)-x|\>|x-a|\wedge|x-b|)\<C\frac{r^{\frac12}}{|x-a|\wedge|x-b|}.
\de
Hence
\be\label{besov1dim4}
\mE[|I_2|^p]\<Cr^{\frac12}.
\ee
In the same way,
\be\label{besov1dim5}
\mE[|I_3|^p]\<Cr^{\frac12}.
\ee
Combining (\ref{besov1dim1})-(\ref{besov1dim5}) gives (\ref{besov1dim0}) and thus finishes the proof.
\end{proof}

Now we can state the Besov regularity of $\s^0_x$.

\bt\label{thm62}
For every $p\>1$, one has $\s^0_x\in \mB_p^\frac{1}{2p}$.
\et
\begin{proof}
Let $B$ and $B'$ be two independent one-dimensional Brownian motions.
Define for $0<u<v$
$$
B^{u,v}_t=\begin{cases}
B_t, &~t\< u,\\
B_u+\theta_uB'_{t-u}, &~u<t\< v,\\
B^{u,v}_v+\theta_vB_{t-v},&~t>v,
\end{cases}
$$
and set
$$
W:=\theta_uB,~~~~~~\hW_t:=\theta_uB^{u,v}.
$$
Then $B^{u,v}$, $W$ and $\hW$ are all Brownian motions.
Hence we can set
$$
X(t,x)=F(x, B)(t),
$$
$$
Y(t,x):=F(x, B^{u,v})(t),
$$
$$
Z(t,x)=F(x, W)(t),
$$
$$
\hZ(t,x)=F(x, \hW)(t).
$$
Define
$$
\a_x:=\inf\{t  \> 0: X(t,x)\in \{a,b\}\},
$$
$$
\b_x:=\inf\{t \> 0: Y(t,x)\in \{a,b\}\}.
$$
Then
$$
\a_x(\om)=\s^0_x(B(\om))~~~\b_x(\om)=\s^0_x(B^{u,v}(\om)).
$$
Hence to prove the theorem it suffices to prove
\be
\mE[|\a_x-\beta_x|^p]\<C(v-u)^{\frac12}.
\ee
Set  $r=v-u$, and
$$
\lambda_x:=\inf\{t \> 0: Z(t,x)\in \{a,b\}\},
$$
$$
\hat{\lambda}_x:=\inf\{t\> 0: \hZ(t,x)\in \{a,b\}\}.
$$
Then
$$
\lambda_x(\om)=\s^0_x(W(\om)),~~~~~\hat{\lambda}_x(\om)=\s^0_x(\hW(\om)).
$$
Since $X(t,x)=Y(t,x)$ for $t\<u$, we have
$$
\a_x1_{\{\a_x\<u\}}=\b_x1_{\{\b_x\<u\}}.
$$
Consequently  $\{\a_x>u\}=\{\b_x>u\}$ and on this set we have $X(u,x)=Y(u,x)\in (a,b)$. Thus
$$
\a_x=u+\lambda_{X(u,x)},~~~\b_x=u+\hat{\lambda}_{Y(u,x)}.
$$
Consequently we have by the strong Markov property and Lemma \ref{lemma61} that
\ce
&&\mE[|\a_x-\b_x|^p]\\
&=&\mE[|\a_x-\b_x|^p1_{\{\a_x>u\}}]\\
&=&\mE[|\lambda_{X(u,x)}-\hat{\lambda}_{Y(u,x)}|^p1_{\{\a_x>u\}}]\\
&=&\mE[\mE[|\lambda_y-\hat{\lambda}_y|^p]_{y=X(u,x)}1_{\{\a_x>u\}}]\\
&\<&Cr^\frac12
\de
as desired.
\end{proof}

Next we turn to the truncated exit time in the multi-dimensional case and suppose that we are in the situation of Section 4. Let $F$ be the unique strong solution to the Equation $e(\s,b)$. Let $T>0$ be arbitrarily fixed and set
$$
\s^0_x(w):=\inf\{t\> 0: F(x,w)(t)\in\partial D\}\wedge T.
$$
Let $B$ and $B'$ be two independent $m$-dimensional Brownian motions and
define
$$
\hB_t=\begin{cases}
B'_t, &~0\< t\< r,\\
B'_r+\theta_rB_{t-r},&~t>r,
\end{cases}
$$
and set
$$
X(t,x,\om):=F(x, B(\om))(t),
$$
$$
\hX(t,x,\om):=F(x, \hB(\om))(t).
$$
Then $X$ and $\hX$ satisfy respectively
\be\label{eq502}
\begin{cases}
&dX(t)=\s_i(X(t))dB^i(t)+b(X(t))dt,\\
&X(0)=x,
\end{cases}
\ee
and
\be\label{eq503}
\begin{cases}
&d\hX(t)=\s_i(\hX(t))d\hB^i(t)+b(\hX(t))dt,\\
&\hX(0)=x,
\end{cases}
\ee
Define
$$
\s_x(\om):=\inf\{t\> 0: X(t,x,\om)\in\partial D\}\wedge T,
$$
$$
\hs_x(\om):=\inf\{t\> 0: \hX(t,x,\om)\in\partial D\}\wedge T,
$$
Then
$$
\s_x(\om)=\s^0_x(B(\om)),~~~~\hs_x(\om)=\s^0_x(\hB(\om)).
$$

We first state the following  counterpart of Lemma \ref{lemma61}.

\bl
For every $p\>1$ there exists a constant $C$ such that for all $x\in \overline{D}$ and $r\<1$
$$
\mE[|\s_x-\hs_x|^p]\<Cr^\frac12.
$$
\el

\begin{proof}
We give a proof for completeness although it is similar to the proof of Lemma \ref{lemma61}. Set
$$
\tB:=\theta_rB,
$$
$$
\tX(t,x,\om):=F(x, \tB(\om))(t),
$$
and
$$
\ts_x:=\inf\{t\> 0: \tX(t,x,\om)\in\partial D\}\st{\wedge T}.
$$
Then we have
\be\label{eq504}
\begin{cases}
&d\tX(t)=\s_i(\tX(t))d\tB^i(t)+b(\tX(t))dt,\\
&\tX(0)=x,
\end{cases}
\ee
and
$$
\ts_x(\om)=\s_x^0(\tB(\om)).
$$
Now observe that
\be\label{l530}
|\s_x-\hs_x|&\<&|\ts_{X(r,x)}-\ts_{\hX(r,x)}|1_{\{\s_x\wedge\hs_x>r\}}+|\s_x|1_{\{\hs_x\<r,\s_x>r\}}\no\\
&&+|\hs_x|1_{\{\s_x\<r,\hs_x>r\}}+r1_{\{\s_x\vee\hs_x\<r\}}\no\\
&=:&\sum_{i=1}^4 I_i.
\ee
By the strong Markov property and Theorem \ref{main result multidimensional truncated} we have
\be\label{l531}
\mE[I_1^p]&\<& \mE[\mE[|\ts_y-\ts_z|^p]|_{y=X(r,x), z=\hX(r,x)}]\no\\
&\<& C\mE[|y-z||_{y=X(r,x), z=\hX(r,x)}]\no\\
&\<&C\mE[|X(r,x)-\hX(r,x)|]\no\\
&\<&Cr^\frac12.
\ee
It is obvious that
\be\label{l532}
\mE[I_4^p]\<Cr^p.
\ee
It remains to estimate $I_2$ and $I_3$. We have by the independence of $B$ and $B'$
\ce
\mE[I_2^p]&=&\mE[|\s_x|^p1_{\{\hs_x\<r,\s_x>r\}}]\no\\
&=&\mE[|\s_x|^p1_{\{\s_x>r\}}]\mE[1_{\{\hs_x\<r\}}]\no\\
&\<&\mE[|\s_x|^p]P(\hs_x\<r).
\de
By Corollary \ref{statement:general_upper_bound_EtauxTp} we have
\ce
\mE[|\s_x|^p]\<C d(x,\partial D),
\de
while
\ce
\mP(\hs_x\<r)\<\mP(\max_{0\<s\<r}|\hX(s,x)-x|\>d(x,\partial D))\<C\frac{r^\frac12}{d(x,\partial D)}.
\de
Hence
\be\label{l533}
 \mE[I_2^p]\<Cr^\frac12.
\ee
In the same way
\be\label{l534}
 \mE[I_3^p]\<Cr^\frac12.
\ee
Combining (\ref{l530})-(\ref{l534}) completes the proof.
\end{proof}

Using this lemma and the same argument as in the proof of Theorem \ref{thm62}, we can prove

\bt
For $p\>1$ one has
$$
\s_x^0\in \mB_p^\frac{1}{2p}.
$$
\et

The regularity obtained is optimal, as the following example illustrates.

Let $w$ be a Brownian motion, and let $x\in [-1,1]$. Define
$$
\s_x:=\s_x(w):=\inf\{t \> 0: |x+w(t)|=1\}.
$$

We first have the following

\bl\label{ll}
For every $p>1$, there exist constants $c,\delta>0$ such that for $|x-y|\< \delta$ it holds
$$
\mE[|\s_x-\s_y|^p]\> c|x-y|.
$$
\el
\begin{proof}
Suppose without any loss of generality that $y<x$. Take $\delta$ as appearing in Lemma \ref{auxi}. Then
for $|x-y|<\delta$ we have
\ce
|\s_x-\s_y|&=&\s_x\vee\s_y-\s_x\wedge\s_y\\
&=&\s_{(1-x+y)1_{\{\s_x<\s_y\}}+(-1+x-y)1_{\{\s_y<\s_x\}}}(\theta_{\s_x\wedge\s_y}w).
\de
Noticing that $\{\s_x<\s_y\}, \{\s_y<\s_x\}\in\cF_{\s_x\wedge\s_y}$, we have by the strong Markov property
and Lemma \ref{auxi}
\ce
\mE[|\s_x-\s_y|^p]&=&\mE[\mE[|\s_z|^p]|_{z=(1-x+y)1_{\{\s_x<\s_y\}}
+(-1+x-y)1_{\{\s_x>\s_y\}}}]\\
&\>&c\mE[\mE[|z-1|\wedge|z+1|]|_{z=(1-x+y)1_{\{\s_x<\s_y\}}
+(-1+x-y)1_{\{\s_x>\s_y\}}}]\\
&\>& c\mE[|x-y|]\\
&=&c|x-y|.
\de
\end{proof}

Let $B$ and $B'$ be two independent $\{\cF_t\}$-Brownian motions defined on some filtered probability space $(\Om,\cF, \{\cF_t\}, \mP)$
satisfying the usual conditions.
For $r>0$ define
$$
\hB_t=\begin{cases}
B'_t, &~t\< r,\\
B'_r+\theta_rB_{t-r},&~t>r.
\end{cases}
$$
Denote
$$
\s:=\inf\{t \> 0: |B_t|=1\}, ~~\hs:=\inf\{t \> 0: |\hB_t|=1\}.
$$
We have:
\bp
There exists $c>0$ such that for $r$ sufficiently small
$$
\mE[|\s-\hs|^p]\> cr^{\frac12}.
$$
\ep
\begin{proof}
Set
$$
\ts_x:=\inf\{t \>0: |x+(\theta_rB)_t|=1\}.
$$
We now take the same $\delta$ as in  Lemma \ref{ll}.
Then we have by Lemma \ref{ll} that
\ce
&&\mE[|\s-\hs|^p]\\
&\> &\mE[|\s-\hs|^p 1_{\{\s\wedge\hs>r\}}1_{\{|B_r-B'_r|\<\delta\}}]\\
&\>&\mE[\mE[|\ts_{B_r}-\ts_{B'_r}|^p1_{\{\s\wedge\hs>r\}}1_{\{|B_r-B'_r|\<\delta\}}|\cF_r]]\\
&=&\mE[1_{\{\s\wedge\hs>r\}}1_{\{|B_r-B'_r|\<\delta\}}\mE[|\ts_x-\ts_y|^p]|_{x=B_r, y=B'_r}]\\
&\>& c\mE[1_{\{\s\wedge\hs>r\}}1_{\{|B_r-B'_r|\<\delta\}}|B_r-B'_r|]\\
&\>& c\{\mE[1_{\{|B_r-B'_r|\<\delta\}}|B_r-B'_r|]-\mE[1_{\{\s\wedge\hs\< r\}}1_{\{|B_r-B'_r|\<\delta\}}|B_r-B'_r|]\}\\
&\>& c_1r^{\frac12}-c_2\mE[|B_r-B'_r|^2]^\frac12\mP(\s\wedge\hs\< r)^\frac12\\
&\>& c_1r^\frac12-c_2r^\frac12r^\frac12\\
&=&c_1r^\frac12-c_2r\\
&\>&cr^\frac12
\de
for $r $ sufficiently small.
\end{proof}

The above estimate obviously implies
\bt
For  $p\>1$ and $\a>\frac{1}{2p}$, $\s\notin \mB_p^{\a}$.
\et

{\bf Acknowledgement}. It is a pleasure to thank Stefan Geiss and Jicheng Liu for fruitful discussions.
\vskip 0.3cm
{\bf Disclaimer}. The authors did not use any AI-tools, neither to obtain the results, nor to complete the manuscript.

\end{document}